\documentstyle[11pt]{article}
\newcommand{\too}{\longrightarrow}
\newcommand{\om}{\omega}

\newcommand{\Li}{{\cal L}}

\newcommand{\G}{{\cal G}}

\newcommand{\Om}{\Omega}
\newcommand{\na}{\nabla}
\newcommand{\wi}{\widetilde}
\newcommand{\al}{\alpha}
\newcommand{\be}{\beta}

\newcommand{\Ga}{\Gamma}

\def \reel{ {\rm I}\!{\rm R} }

 \def \rat{ {\rm Q}\kern-.65em {}^{{}_/ }}

 \def\ent{{{\rm Z}\mkern-5.5mu{\rm Z}}}

\newtheorem{th}{Theorem}[section]
\newtheorem{pr}{Proposition}[section]
\newtheorem{Le}{Lemma}[section]

\newtheorem{exem}{Example}[section]

\title{Polynomial Poisson structures on  affine solvmanifolds}
\author{Mohamed Boucetta- Alberto Medina\footnote{Recherche men\'ee dans le
cadre du projet de coop\'eration CNRS-CNRST " Vari\'et\'e de
Poisson et g\'eom\'etrie affine" SPM07/05.}} \date{
}\parindent=0cm
\begin{document}
\maketitle

{\bf Abstract.} A $n$-dimensional Lie group $G$ equipped with a
left invariant symplectic form $\om^+$ is called a symplectic Lie
group. It is well-known that $\om^+$ induces a left invariant
affine structure on $G$. Relatively to this affine structure we
show that the left invariant Poisson tensor $\pi^+$ corresponding
to $\om^+$ is polynomial of degree 1 and any right invariant
$k$-multivector field on $G$ is polynomial of degree at most $k$.
If $G$ is unimodular, the symplectic form $\om^+$ is also
polynomial and the volume form $\wedge^{\frac{n}2}\om^+$ is
parallel. We show also that any left invariant  tensor field on a
nilpotent symplectic Lie group is polynomial, in particular, any
left invariant Poisson structure on a nilpotent symplectic Lie
group is polynomial. Because many symplectic Lie groups admit
uniform lattices, we get a large class of polynomial Poisson
structures on compact affine solvmanifolds.\bigskip



 \section{Introduction and main results}

Recall that an affine manifold is a differential manifold $M$
together with a special atlas of coordinate charts such that all
coordinate changes extend to affine automorphisms of $\reel^n$.
These distinguished charts are called affine charts. The data of a
flat and  torsion free connection $\na$ on a manifold $M$ is
equivalent to the data of an affine structure.

A tensor field on an affine manifold $M$ is called {\it
polynomial} if in  affine coordinates its coefficients are
polynomial functions.  A Poisson structure on an affine manifold
is called polynomial if the space of local polynomial functions is
closed under the Poisson bracket. In an equivalent way this means
that the associated Poisson bivector is polynomial. For some
general results on polynomial tensor fields see [4,6,7,8,16]. Let
us describe briefly the affine structure associated to a Lie group
endowed with a left invariant symplectic form. This affine
structure is the context on which we will  state our main results
on the polynomial nature of some tensor fields and some Poisson
structures.

 Let $G$ be a Lie
group with Lie algebra $\G=T_eG$, where $e$ stands for the unit
 of $G$. For any tensor $T$ on $\G$, we denote by $T^+$ and $T^-$
 respectively the left invariant tensor field and the right
 invariant tensor field on $G$ associated to $T$.  If $\om$ is a scalar
  non degenerate 2-cocycle of
 $\G$, the differential 2-form $\om^+$ on $G$
 is a left invariant symplectic form on $G$ and $(G,\om^+)$ is
 called a symplectic Lie group. Symplectic Lie groups were studied by several authors see,
  for instance, [1,2,3,10,11,13,14]. A connected Lie group $G$ is
  symplectic if and only if its universal covering $\hat G$ admits
  an etale  representation by affine transformations of $\G^*$
  with linear part the coadjoint representation of $\hat G$ and
  infinitesimal part a skew-symmetric 1-cocycle (see [12]). This
  implies that
  the formula
\begin{equation}\om^+(\na_{u^+}v^+,w^+)=-\om^+(v^+,[u^+,w^+]),
\end{equation}where $u,v,w\in\G$,  defines a left invariant flat and torsion free connection
$\na$. This affine structure will be called the affine structure
associated to the symplectic Lie group $(G,\om^+)$.\bigskip

Let us state our mains results.
\begin{th} Let $(G,\om^+)$ be a  connected symplectic Lie group of dimension
$n$ endowed with the associated affine structure. Then the
following assertions hold.
\begin{enumerate}\item In
a neighborhood of any element of $G$, there exists an affine chart
$(x_1,\ldots,x_n)$ such that, for any $i,j=1,\ldots,n$, the
Poisson bracket of $x_i$ and $x_j$ associated to $\om^+$ is given
by
$$\{x_i,x_j\}=\sum_{k=1}^nC_{ij}^kx_k+\mu_{ij},$$where $C_{ij}^k$
are constants of structure of the Lie algebra of $G$ and
$\mu_{ij}$ are  constants.\item Any right invariant
$k$-multivector field on $G$ is polynomial of degree at most $k$.
\item If $G$ is unimodular then the symplectic form $\om^+$ is
polynomial of degree at most $n-1$, the volume form
$\wedge^{\frac{n}2}\om^+$ is parallel and any right invariant
differential form on $G$ is polynomial.
\end{enumerate}
\end{th}

\begin{th} Let $(G,\om^+)$ be a connected nilpotent  symplectic Lie group
endowed with the associated affine structure. Then any left
invariant multivector field on $G$ is polynomial. In particular,
any left invariant Poisson structure on $G$ is polynomial.\end{th}

There are some interesting implications of Theorems 1.1 and 1.2.
\begin{enumerate} \item Let $(G,\om^+)$ be a connected
$n$-dimensional symplectic Lie group. If $G$ admits an uniform
 lattice, it is well known that $G$ is unimodular. On the other hand,
 according to a result of Medina-Lichnerowicz [11], the associated
 affine structure to $(G,\om^+)$ is geodesically complete if and
 only if $G$ is unimodular and, in this case $G$ is solvable.

Consequently if $\Ga$ is an (uniform) lattice in $G$ then
$M={\Ga}\setminus G$ is a
 compact solvmanifold
  which carries an affine structure and a symplectic form
 $\wi\om$ such that:
 \begin{enumerate}\item the Poisson bracket corresponding to $\wi\om$
 is polynomial of degree 1;
 \item the symplectic form $\wi\om$ is polynomial of degree at
 most $n-1$ and the volume form $\wedge^{\frac{n}2}\wi\om$ is
parallel. Note that a compact affine manifold with a parallel
volume form possesses interesting properties (see
[9]).\end{enumerate}\item Let $(G,\om^+)$ be a connected nilpotent
symplectic Lie group. Then: \begin{enumerate}\item any exact
Lie-Poisson tensor on $G$ is polynomial;\item for any  uniform
 lattice $\Ga$  on $G$, any solution of the classical Yang-Baxter equation on
 $\G$ gives
 arise to a left
invariant Poisson tensor on $G$ which
 projects on a polynomial Poisson tensor  on $M=\Ga\setminus
 G$.\end{enumerate}
 \item Let $(G,\om^+)$ be a symplectic Lie group
and let $r\in\G\wedge\G$ be the solution of the classical
Yang-Baxter equation associated to $\om$. According to Theorem
1.1, $r^+$ is a polynomial Poisson structure of degree 1 and $r^-$
is a polynomial Poisson structure of degree at most 2. Thus, we
recover  a result of Diatta-Medina (see [5]) which states that the
Lie-Poisson  bivector $r^+-r^-$ is polynomial of degree 2.
\end{enumerate}

The paper is organized as follows. In Section 2, we give some
properties of the affine structure associated to a symplectic Lie
group. Proofs of Theorems 1.1 and 1.2 are developed in Section 3.
Using Theorem 1.1 and results of Medina-Revoy on lattices in
symplectic Lie groups [14], we exhibit in Section 4 and infinity
of non homeomorphic compact affine solvmanifolds endowed with
polynomial Poisson tensors.

\bigskip

{\bf Acknowledgment }  This work was finalized during the stay of
the first author at the university Montpelier II. The first author
thanks the department of mathematics for having invited him.

\section{Some properties of the affine structure associated to a
symplectic Lie group}

This section  is a preparation of Section 3 in which we will prove
Theorems 1.1 and 1.2.

First, we will consider the affine structure given by $(1)$ from a
different point of view which will be useful through the paper.

Let $\pi^+$ be a left invariant Poisson bivector on a Lie group
$G$ and  denote by $\pi^+_\#:T^*G\too TG$ the associated
homomorphism. Recall that the Koszul bracket associated to $\pi^+$
is given by
$$[\al,\be]_{\pi^+}={\cal L}_{\pi^+_\#(\al)}\be-{\cal
L}_{\pi^+_\#(\be)}\al-d\pi^+(\al,\be),$$ where $\al$ and $\be$ are
differential 1-forms on $G$ and ${\cal L}$ denotes the Lie
derivative. This bracket endows $\Om^1(G)$ with a structure of a
Lie algebra and, for any $\al,\be\in\Om^1(G)$,
$$
{\pi^+_\#}\left([\al,\be]_{\pi^+}\right)=[{\pi^+_\#}(\al),{\pi^+_\#}(\be)].$$

An easy calculation gives that, for any vector field $X$ on $G$
and for any differential 1-form $\al$,
\begin{equation} [{\pi^+_\#}(\al),X]=-{\cal
L}_X\pi^+(\al,.\;)-{\pi^+_\#}({\cal L}_{X}\al).\end{equation} We
deduce easily from $(2)$ that for any right invariant 1-forms
$\al^-$ and $\be^-$ and for any right invariant vector field
$X^-$, $$ [\al^-,\be^-]_{\pi^+}(X^-)=-{\cal
L}_{X^-}\pi^+(\al^-,\be^-)=0.$$ Thus, for  any right invariant
1-forms $\al^-$ and $\be^-$, the bracket $[\al^-,\be^-]_{\pi^+}$
vanishes and hence
\begin{equation}[{\pi^+_\#}(\al^-),{\pi^+_\#}(\be^-)]=0.\end{equation}

With this remark in mind, we consider   a connected symplectic Lie
group $(G,\om^+)$  and we denote by $\pi^+$ the associated left
invariant Poisson tensor. From $(3)$, for any basis
$(\al_1^-,\ldots,\al_n^-)$ of right invariant 1-forms on $G$,
$(\pi^+_\#(\al^-_1),\ldots,\pi^+_\#(\al^-_n))$ is a commuting
parallelism
 of vector fields on
$G$. This  defines a flat and torsion free linear connection
$\wi\na$ by putting
$$\wi\na\pi^+_\#(\al^-)=0$$for any right invariant
1-form on $G$. Moreover, for any right invariant vector field
$X^-$, we get from $(2)$ $$ [{\pi^+_\#}(\al^-),X^-]=
\wi\na_{{\pi^+_\#}(\al^-)}X^-=-{\pi^+_\#}({\cal
L}_{X^-}\al^-),$$and since ${\cal L}_{X^-}\al^-$ is a right
invariant 1-form, we get
$$\wi\na\wi\na_{{\pi^+_\#}(\al^-)}X^-=0.$$Then $X^-$ is
an affine infinitesimal transformation and hence $\wi\na$ is a
left invariant linear connection.

Now, let us compare the linear connection $\na$ defined by $(1)$
and $\wi\na$. More precisely, we will show that
$$\wi\na=\na.$$

The symplectic form $\om^+$ gives arise to an  isomorphism
$\om^\flat:TG\too T^*G$, $u\mapsto \om^+(u,.)$. This isomorphism
and its inverse $(\om^\flat)^{-1}:T^*G\too TG$ define an
isomorphism between the space of tensor field of type $(p,q)$ and
the space of tensor field of type $(q,p)$. For any tensor field
$T$, we denote by $T^\om$ its image under this isomorphism. For
instance, for any vector field $X$ and any 1-form $\al$, $X^\om$
is the 1-form $i_X\om^+$ and $\al^\om=-\pi^+_\#(\al)$. It is
obvious that $(T^\om)^\om=T.$
\begin{pr} Let $(G,\om^+)$ be a symplectic Lie group. Then:
\begin{enumerate}\item For any left invariant vector field $X^+$ and for any tensor field
$T$, we have $$\na_{X^+}T=\left( {\cal L}_{X^+}
T^\om\right)^\om,$$where ${\cal L}_{X^+}$ is the Lie derivative in
the direction of $X^+$.

 \item A tensor field
$T$ is parallel with respect to $\na$ if and only if $T^\om$ is
right invariant.\end{enumerate}
\end{pr}

{\bf Proof.} Note that the second assertion  is an immediate
consequence of the first one. Let us establish the first
assertion. Since, for any tensor fields $T_1$ and $T_2$,
$(T_1\otimes T_2)^\om=T_1^\om\otimes T_2^\om$ and since both
${\cal L}_{X^+}$ and $\na_{X^+}$ are derivative, it suffices to
establish the relation for left invariant vector fields and left
invariant differential 1-forms. Let $Y^+$ and $Z^+$ be  left
invariant vector fields. We have
\begin{eqnarray*} \om^+(\left( {\cal L}_{X^+}
i_{Y^+}\om^+\right)^\om,Z^+)&=&{\cal L}_{X^+} i_{Y^+}\om^+(Z^+)\\
&=&X^+.\om^+(Y^+,Z^+)-\om^+(Y^+,[X^+,Z^+])\\
&=&-\om^+(Y^+,[X^+,Z^+])=\om^+(\na_{X^+}Y^+,Z^+),\end{eqnarray*}and
the formula follows. One can deduce easily the formula for a left
invariant 1-form.$\Box$\bigskip

An immediate consequence of this proposition is that,
$\na{\pi^+_\#}(\al^-)=0$  for any right invariant 1-form $\al^-$
and hence $\wi\na=\na$.
\bigskip

Now, given a connected symplectic Lie group $(G,\om^+)$,  let us
construct an affine   atlas corresponding $\na$.

Since $\om^+$ is left invariant , for any $u\in \G$, the vector
field $u^-$ is symplectic, i.e., $$ {\cal
L}_{u^-}\om^+=di_{u^-}\om^+=0.$$ Thus, for any basis
$(u_1,\ldots,u_n)$ of $\G$, there exists in a neighborhood of any
element of $G$ a local coordinates  $(x_1,\ldots,x_n)$  such that
\begin{equation} i_{u_i^-}\om^+=dx_i,\quad
i=1,\ldots,n.\end{equation}
 We get from Proposition 2.1
that $\na dx_i=0$ and  we deduce that $(x_1,\ldots,x_n)$ are
affine coordinates.

Now we will express  $\om^+$ and  $\pi^+$ in the  affine
coordinates constructed above. Fix a basis $(u_1,\ldots,u_n)$ of
$\G$, denote by $(\al_1,\ldots,\al_n)$ its dual basis and consider
the affine coordinates $(x_1,\ldots,x_n)$ given by $(4)$.  In this
coordinates we have
\begin{eqnarray}
\om^+&=&\sum_{i<j}\pi^+(\al_i^-,\al_j^-)dx_i\wedge dx_j,\\
\pi^+&=&\sum_{i<j}\om^+(u_i^-,u_j^-)\frac{\partial}{\partial
x_i}\wedge \frac{\partial}{\partial x_j},\\
 u^-&=&\sum_{j=1}^n\om^+(u_j^-,u^-)\frac{\partial}{\partial
x_j},\quad u\in\G,\\
\al^-&=&\sum_{j=1}^n\pi^+(\al_j^-,\al^-)dx_j,\quad\al\in\G^*.\end{eqnarray}

The following proposition will play a crucial role in the proof of
Theorem 1.1.
\begin{pr} Let $(G,\om^+)$ be a connected symplectic Lie group
endowed with the associated affine structure. Then, for any
$u,v\in\G$ such that $[u,v]\not=0$ , $\om^+(u^-,v^-)$ is a
polynomial function of degree  1.\end{pr}

{\bf Proof.} The right invariant vector fields $u^-,v^-$ are
symplectic and hence $[u^-,v^-]$ is hamiltonian and we have
$$i_{[u^-,v^-]}\om^+=d\om^+(u^-,v^-).$$ By using
Proposition 2.1, we get $ \na d\om^+(u^-,v^-)=0$ and the result
follows.$\Box$\bigskip

Note that if $A$ is the matrix
$\left(\om^+(u_i^-,u_j^-)\right)_{1\leq i,j\leq n}$, we have
\begin{equation} A^{-1}=\left(\pi^+(\al_i^-,\al_j^-)\right)_{1\leq i,j\leq
n}.\end{equation}

 Consider now the volume form
$\Om^+=\wedge^{\frac{n}2}\om^+.$ We have, from $(7)$,
\begin{equation} \Om^+(u_1^-,\ldots,u_n^-)=\det(A)\Om^+(\frac{\partial}{\partial
x_1},\ldots,\frac{\partial}{\partial x_n}),\end{equation} and
\begin{equation}\det(A)=\frac1{\left((\frac{n}2)!\right)^2}\Om^+(u_1^-,\ldots,u_n^-)^2.
\end{equation}

\begin{pr} Let $(G,\om^+)$ be a connected unimodular  symplectic Lie group
endowed with the associated affine structure and let $\pi^+$ be
the left invariant Poisson structure corresponding to $\om^+$.
Then $\pi^+(\al^-,\be^-)$ is a polynomial function of degree at
most $n-1$, for any  differential forms $\al^-,\be^-$.\end{pr}

{\bf Proof.} Since $G$ is  unimodular, $\Om^+$ is a right
invariant and then $\Om^+(u_1^-,\ldots,u_n^-)$ is constant. Hence,
from $(11)$, $\det(A)$ is a constant. On the other hand, from
Proposition 2.2, the coefficients of $A$ are polynomial functions
of degree 1, consequently the coefficients of the inverse $A^{-1}$
are polynomial of degree at most $n-1$ and the proposition follows
from $(9)$.$\Box$\bigskip

 The following Lemma will be useful in the proof of Theorem
1.2.
\begin{Le} A function $f$ on an  unimodular symplectic Lie group
$(G,\om^+)$ is polynomial if and only if $u^-(f)$ is polynomial
for any $u\in\G$.\end{Le}

{\bf Proof.} We have from $(7)$
\begin{equation}u^-(f)=\sum_{j=1}^n\om^+(u_j^-,u^-)\frac{\partial f}{\partial
x_j}.\end{equation}So if $f$ is polynomial, $u^-(f)$ is polynomial
according to Proposition 2.2. For the converse, we deduce from
$(12)$ that
$$\left(\begin{array}{c}\frac{\partial f}{\partial
x_1}\\\vdots\\\frac{\partial f}{\partial x_n}\end{array}\right)=
A^{-1}\left(\begin{array}{c}u^-_1(f)\\\vdots\\u_n^-(f)\end{array}\right).$$

But, we have see in the proof of Proposition 2.3 that if $G$ is
unimodular the coefficients of $A^{-1}$ are polynomial and the
Lemma follows.$\Box$

\section{Proof of Theorems 1.1 and  1.2}

{\bf Proof of Theorem 1.1.} \begin{enumerate}\item Fix a basis
$(u_1,\ldots,u_n)$ of $\G$, denote by $(\al_1,\ldots,\al_n)$ its
dual basis and consider the affine  coordinates $(x_1,\ldots,x_n)$
given by $(4)$. For any $u_i,u_j$, put
$[u_i,u_j]=\sum_{k=1}^nC_{ij}^ku_k$.

 From $(6)$,  we have
$\{x_i,x_j\}=\om^+(u_i^-,u_j^-)$. Hence, for any $1\leq i,j,k\leq
n$,
\begin{eqnarray*}
\frac{\partial}{\partial
x_k}.\{x_i,x_j\}&=&d\om^+(u_i^-,u_j^-)(\pi_{\#}^+(\al_k^-))\\
&=&i_{[u_i^-,u_j^-]}\om^+(\pi_{\#}^+(\al_k^-))\\
&=&\al_k^-([u_i^-,u_j^-])=C_{ij}^k,
\end{eqnarray*}thus $d\{x_i,x_j\}=\sum_{k=1}^nC_{ij}^kdx_k$ and
the desired relation follows.

\item From $(7)$ and Proposition 2.2 we deduce that any right
invariant vector field on $G$ is polynomial of degree at most 1
and hence any right invariant $k$-multivector field must be
polynomial of degree at most $k$. \item This is an immediate
consequence of Proposition 2.3, $(5)$, $(9)$ and $(10)$. $\Box$

\end{enumerate}

{\bf Proof of Theorem 1.2}  To prove the theorem, it suffices to
show that if $G$ is nilpotent then any left invariant vector field
is polynomial. Fix a basis $(u_1,\ldots,u_n)$ of $\G$ and consider
the affine  coordinates $(x_1,\ldots,x_n)$ given by $(4)$. Let
$u^+$ be a left invariant vector field on $G$. We have
$$u^+=\sum_{i=1}^n\om^+(u_i^-,u^+)\frac{\partial}{\partial x_i}.$$

Since $\Li_{u_j^-}\om^+=0$ and $[u_j^-,u^+]=0$, we have, for any
$1\leq i,j\leq n$,
$$u^-_j.\om^+(u_i^-,u^+)=\om^+([u_j^-,u_i^-],u^+)$$and by
induction, we get
$$u_{j_1}^-\circ\ldots\circ u_{j_r}^-.\om^+(u_i^-,u^+)=
\om^+([u^-_{j_1},[\ldots,[u_{j_r}^-,u_i^-]],u^+)$$for any $1\leq
j_1,\ldots,j_r\leq n$. Since $G$ is nilpotent, we get for $r$
large$$u_{j_1}^-\circ\ldots\circ u_{j_r}^-.\om^+(u_i^-,u^+)=0$$
and we deduce from Lemma 2.1 that $\om^+(u_i^-,u^+)$ is a
polynomial function and the theorem is proved.$\Box$

\section{Examples}

In this section, we give a large class of four dimensional
solvmanifolds which admit polynomial symplectic forms and
polynomial Poisson structures. The construction is based on
Theorems 1.1 and 1.2 and on the results of Medina-Revoy on
Lattices in four dimensional symplectic Lie groups (see [14,15]).
Indeed, in [12] Medina and Revoy showed that there are four non
abelian real unimodular Lie algebras of dimension four endowed
with a scalar non degenerate 2-cocycle. Moreover, the connected
and simply connected Lie group of any of such Lie algebras has an
infinity of non isomorphic lattices. For any Lie algebra in the
list of Medina-Revoy, we consider the corresponding connected and
simply connected Lie group endowed with  a left invariant
symplectic form,  we give a global affine  chart and we express
the symplectic form and the Poisson bivector in this chart.
Finally, we give a description of lattices in this group.

\begin{exem}\begin{enumerate}\item

We consider the  Lie group $G_1=\reel^4$ with the product
$$(x,y,z,t)(x',z',y',t')=(x+x',y+y'+tx',z+z'+ty'+\frac{t^2}2x',t+t').$$

We denote by $\G_1$ its Lie algebra and by $(e_1,e_2,e_3,e_4)$ the
canonical basis of $\G_1$. We have
$$e^+_1=\partial_x+t\partial_y+\frac{t^2}2\partial_z,\;
e_2^+=\partial_y+t\partial_z,\; e_3^+=\partial_z,\;
e_4^+=\partial_t.$$
$$e_1^-=\partial_x,\; e_2^-=\partial_y,\; e_3^-=\partial_z,\;
e_4^-=x\partial_y+y\partial_z+\partial_t.$$
 The nonzero brackets of $\G_1$ are the following:
$$[e_4,e_1]=e_2\quad\mbox{and}\quad [e_4,e_2]=e_3.$$
We consider the scalar non degenerate 2-cocycle on $\G_1$ given by
$\om=e_4^*\wedge e_3^*+e_1^*\wedge e_2^*.$ A direct computation
gives that the corresponding symplectic 2-form on $G_1$ is given
by
$$\om^+=dx\wedge dy-\frac12t^2dx\wedge dt+tdy\wedge dt-dz\wedge
dt,$$and
\begin{eqnarray*} i_{e_1^-}\om^+&=&d(y-\frac16t^3),\quad
i_{e_2^-}\om^+=d(-x+\frac12t^2),\\
i_{e_3^-}\om^+&=&-dt,\quad
i_{e_4^-}\om^+=d(-\frac12x^2+\frac12xt^2-yt+z).\end{eqnarray*} We
put
$$X=y-\frac16t^3,\; Y=-x+\frac12t^2,\; Z=-t\quad\mbox{and}\quad T=
-\frac12x^2+\frac12xt^2-yt+z,$$and we get an affine chart
$(X,Y,Z,T)$.

Here the matrix $A=\left(\om^+(e_i^-,e_j^-)\right)_{1\leq i,j\leq
4}$ is given by
$$A=\left(\begin{array}{cccc}0&1&0&-Y\\-1&0&0&-Z\\0&0&0&-1\\Y&Z&1&0\end{array}\right)$$
and hence, by using $(6)$,
$$\pi^+=\partial_X\wedge
\partial_Y-Y\partial_X\wedge\partial_T-Z\partial_Y\wedge\partial_T-\partial_Z\wedge\partial_T.$$
The inverse of $A$ is given by
$$A^{-1}=\left(\begin{array}{cccc}0&-1&Z&0\\1&0&-Y&0\\-Z&Y&0&1\\0&0&-1&0\end{array}\right)$$
and hence
$$\om^+=-dX\wedge dY+ZdX\wedge dZ-YdY\wedge dZ+dZ\wedge dT.$$

It is clear that $\Ga=\{(m,n,k,2r)/ m,n,k,r\in\ent\}$ is an
uniform lattice in $G_1$. Actually, $G_1$ admits an infinity non
isomorphic lattices.

\item We consider the  Lie group $G_2=\reel^4$ with the product
$$(x,y,z,t)(x',y',z',t')=(x+x',y+e^xy',z+e^{-x}z',t+t').$$
We denote by $\G_2$ its Lie algebra and by $(e_1,e_2,e_3,e_4)$ the
canonical basis of $\G_2$. We have
$$e_1^+=\partial_x,\; e_2^+=e^x\partial_y,\;
e_3^+=e^{-x}\partial_z,\; e_4^+=\partial_t,$$
$$e_1^-=\partial_x+y\partial_y-z\partial_z,\; e_2^-=\partial_y,\;
e_3^-=\partial_z,\; e_4^-=\partial_t.$$

The nonzero brackets of $\G_2$ are the following:
$$[e_1,e_2]=e_2\quad\mbox{and}\quad [e_1,e_3]=-e_3.$$
We consider the scalar non degenerate 2-cocycle on $\G_2$ given by
$\om=e_1^*\wedge e_4^*+e_2^*\wedge e_3^*.$ A direct computation
gives that the corresponding symplectic 2-form on $G_2$ is given
by
$$\om^+=dx\wedge dt+dy\wedge dz$$and
$$i_{e_1^-}\om^+=d(t+yz),\; i_{e_2^-}\om^+=dz,\;
i_{e_3^-}\om^+=-dy\;\mbox{and}\; i_{e_4^-}\om^+=-dx.$$We put
$(X,Y,Z,T)=(t+yz,z,-y,-x)$ and we get an affine chart.

Here the matrix $A=\left(\om^+(e_i^-,e_j^-)\right)_{1\leq i,j\leq
4}$ is given by
$$A=\left(\begin{array}{cccc}0&-Y&-Z&1\\Y&0&1&0\\Z&-1&0&0\\-1&0&0&0\end{array}\right)$$
and hence, by using $(6)$,
$$\pi^+=-Y\partial_X\wedge
\partial_Y-Z\partial_X\wedge\partial_Z+\partial_X\wedge\partial_T+\partial_Y\wedge \partial_Z.$$
The inverse of $A$ is given by
$$A^{-1}=\left(\begin{array}{cccc}0&0&0&-1\\0&0&-1&-Z\\0&1&0&Y\\1&Z&-Y&0\end{array}\right)$$
and hence
$$\om^+=-dX\wedge dT-dY\wedge dZ-ZdY\wedge dT+YdZ\wedge dT.$$

The Lie group $G_2$ is a direct product of $G_2'$ with the abelian
group $\reel$ where the multiplication in $G_2'$ is given by
$$(x,y,z)(x',y',z')=(x+x',y+e^xy',z+e^{-x}z').$$ Hence, if $\Ga_1$
is a lattice in $G_2'$ then  $\Ga=\Ga_1\times\ent$ is a lattice in
$G_2$. But, according to [13], a lattice $\Ga_1$ in $G_2'$ is a
semi-direct product of $\ent$ by $\ent^2$ by an action $\phi$ of
$\ent$ on $\ent^2$ with $\phi(1)$ acts as
$A=\left(\begin{array}{cc}0&-1\\1&n\end{array}\right)$ $(n\geq3)$
in a suitable basis of $\Ga_1\cap N$ where $N$ is nilradical of
$G_2'$.

\item We consider the  Lie group $G_3=\reel^4$ with the product
$$(x,y,z,t)(x',y',z',t')=(x+x',y+y'\cos(x)-z'\sin(x),z+y'\sin(x)+z'\cos(x),t+t').$$
We denote by $\G_3$ its Lie algebra and by $(e_1,e_2,e_3,e_4)$ the
canonical basis of $\G_3$. We have
$$e_1^+=\partial_x,\; e_2^+=\cos x\partial_y+\sin x\partial_z,\;
e_3^+=-\sin x\partial_y+\cos x\partial_z,\; e_4^+=\partial_t,$$
$$e_1^-=\partial_x-z\partial_y+y\partial_z,\; e_2^-=\partial_y,\;
e_3^-=\partial_z,\; e_4^-=\partial_t,$$ The nonzero brackets of
$\G_3$ are the following:
$$[e_1,e_2]=e_3\quad\mbox{and}\quad [e_1,e_3]=-e_2.$$
We consider the scalar non degenerate 2-cocycle on $\G_3$ given by
$\om=e_1^*\wedge e_4^*+e_2^*\wedge e_3^*.$ A direct computation
gives that the corresponding symplectic 2-form on $G_2$ is given
by
$$\om^+=dx\wedge dt+dy\wedge dz$$and
$$i_{e_1^-}\om^+=d(t-yz),\; i_{e_2^-}\om^+=dz,\;
i_{e_3^-}\om^+=-dy\;\mbox{and}\; i_{e_4^-}\om^+=-dx.$$We put
$(X,Y,Z,T)=(t-yz,z,-y,-x)$ and we get an affine chart.

The matrix $A=\left(\om^+(e_i^-,e_j^-)\right)_{1\leq i,j\leq 4}$
is given by
$$A=\left(\begin{array}{cccc}0&Z&-Y&1\\-Z&0&1&0\\Y&-1&0&0\\-1&0&0&0\end{array}\right)$$
and hence, by using $(6)$,
$$\pi^+=Z\partial_X\wedge
\partial_Y-Y\partial_X\wedge\partial_Z+\partial_X\wedge\partial_T+\partial_Y\wedge \partial_Z.$$
The inverse of $A$ is given by
$$A^{-1}=\left(\begin{array}{cccc}0&0&0&-1\\0&0&-1&-Y\\0&1&0&-Z\\1&Y&Z&0\end{array}\right)$$
and hence
$$\om^+=-dX\wedge dT-dY\wedge dZ-YdY\wedge dT-ZdZ\wedge dT.$$

The Lie group $G_3$ is a direct product of $G_3'$ with the abelian
group $\reel$ where the multiplication on $G_3'$ is given by
$$(x,y,z)(x',y',z')=(x+x',y+y'\cos(x)-z'\sin(x),z+y'\sin(x)+z'\cos(x)).$$  Hence, if $\Ga_1$
is a lattice in $G_3'$ then $\Ga=\Ga_1\times\ent$ is a lattice in
$G_3$. Now, the Lie group $G_3'$, which  is a semi-direct product
of $\reel$ with $\reel^2$, is isomorphic to the universal covering
of the positive motions of the Euclidian plan. The element
$1\in\reel$ and the standard lattice $\ent^2$ in $\reel^2$
generates a lattice $\Ga_1$ isomorphic to $\ent^3$. Moreover, if
we consider the subgroup of $G_3'$ generated by $\ent^2$ and the
element $\frac12\in\reel$ then we get a lattice in $G_3'$
containing $\Ga_1$ as a subgroup of index 2 which is not
nilpotent.

 \item We consider $G_4=\reel^3\times\reel^{*+}$ the direct
product of the 3-dimensional Heisenberg Lie group with the
1-dimensional Lie group. The product is given by
$$(x,y,z,t)(x',y',z',t')=(x+x',y+y',z+z'+xy',tt').$$
We denote by $\G_4$ its Lie algebra and by $(e_1,e_2,e_3,e_4)$ the
canonical basis of $\G_4$. We have
$$e_1^+=\partial_x,\; e_2^+=\partial_y+x\partial_z,\;
e_3^+=\partial_z,\; e_4^+=t\partial_t,$$
$$e_1^-=\partial_x+y\partial_z,\; e_2^-=\partial_y,\;
e_3^-=\partial_z,\; e_4^-=t\partial_t.$$

The nonzero brackets of $\G_4$ are the following:
$$[e_1,e_2]=e_3.$$
We consider the scalar non degenerate 2-cocycle on $\G_4$ given by
$\om=e_1^*\wedge e_4^*+e_2^*\wedge e_3^*.$ A direct computation
gives that the corresponding symplectic 2-form on $G_4$ is given
by
$$\om^+=\frac1t dx\wedge dt+dy\wedge dz,$$and
$$i_{e_1^-}\om^+=d(\ln t-\frac12y^2),\; i_{e_2^-}\om^+=dz,\;
i_{e_3^-}\om^+=-dy,\; i_{e_4^-}\om^+=-dx.$$ We put $(X,Y,Z,T)=(\ln
t-\frac12y^2,z,-y,-x)$ and we get  an affine chart.

The matrix $A=\left(\om^+(e_i^-,e_j^-)\right)_{1\leq i,j\leq 4}$
is given by
$$A=\left(\begin{array}{cccc}0&Z&0&1\\-Z&0&1&0\\0&-1&0&0\\-1&0&0&0\end{array}\right)$$
and hence, by using $(6)$
$$\pi^+=Z\partial_X\wedge
\partial_Y+\partial_X\wedge\partial_T+\partial_Y\wedge \partial_Z.$$
The inverse of $A$ is given by
$$A^{-1}=\left(\begin{array}{cccc}0&0&0&-1\\0&0&-1&0\\0&1&0&-Z\\1&0&Z&0\end{array}\right)$$
and hence
$$\om^+=-dX\wedge dT-dY\wedge dZ-ZdZ\wedge dT.$$

The group $G_4$ is a product of the Heisenberg group $N_3$ with
the abelian group $\reel$. Note that $N_3$ is isomorphic to the
group of all unipotent real upper triangular $3\times 3$ matrix.
The subset of $N_3$
$$\Ga_{p,q,r}=\left\{\left(\begin{array}{ccc}1&\frac{m}p&\frac{k}{pqr}\\
0&1&\frac{n}q\\0&0&1\end{array}\right);\;m,n,k\in\ent\right\},$$where
$p,q,r$ are fixed integers numbers with $pqr\not=0$, is a lattice
in $N_3$. It is commensurable with $\Ga_{1,1,1}$. In fact,
$\Ga_{p,q,r}$ contains $\Ga_{1,1,1}$ as a subgroup of index
$p^2q^2r$.

\end{enumerate}

\end{exem}

To get examples in dimension greater or equal to 6, one cane use
the results given in [1] or [10]. \bigskip

 {\bf References}\bigskip

[1] {\bf Bajo I.; Benayadi S.; Medina A.,} {\it Symplectic
structures on quadratic Lie algebras,}  J. Algebra {\bf 316}
(2007), no. 1, 174-188.

[2] {\bf Benson C.,  Gordon C.,} {\it K\"ahler and symplectic
structures on nilmanifolds,} Topology {\bf 27, 4} (1988), 513-518.

[3] {\bf Chu B. Y.,} {\it Symplectic homogeneous spaces,} Trans.
Am. Math. Soc. {\bf 197} (1974), 145-159.

[4] {\bf Dekimpe, K.,} {\it Polynomial structures and the
uniqueness of affinely flat infra-nilmanifolds,}  Math. Z. {\bf
224} (1997),  no. 3, 457-481.

[5] {\bf Diatta A. and  Medina A.}, {\it Classical Yang-Baxter
equation and left invariant affine geometry on Lie groups,}
Manuscripta Math. {\bf 114}  (2004),  no. 4, 477-486.

[6] {\bf Fried, D.,} {\it Polynomials on affine manifolds,} Trans.
Amer. Math. Soc. {\bf 274} (1982), no. 2, 709-719.

[7]{\bf Goldman, William M.; Hirsch, Morris W.,} {\it  Polynomial
forms on affine manifolds,}  Pacific J. Math. {\bf  101}  (1982),
no. 1, 115--121.

[8] {\bf Goldman, William M.,} {\it On the polynomial cohomology
of affine manifolds,} Invent. Math. {\bf 65} (1981/82), no. 3,
453-457.

[9] {\bf  Huei Shyong L.,} {\it  Note on affine manifold with
parallel volume element,}  Geom. Dedicata  {\bf 9}  (1980), no. 2,
195-198.

[10] {\bf Khakimdjanov, Yu.; Goze, M.; Medina, A.,} {\it
Symplectic or contact structures on Lie groups,}  Differential
Geom. Appl. {\bf 21} (2004), no. 1, 41-54.

[11] {\bf  Lichnerowicz A. and  Medina A.,} {\it On Lie groups
with left-invariant symplectic or K\"ahlerian structures,}
{Letters in Mathematical Physics {\bf16} (1988), 225-235.}

[12] {\bf  Medina A.,} {\it Structures de Poisson affines,}
Symplectic Geometry and Mathematical Physics, Progress In
Mathematics, Birkh\"auser (1991).

[13] {\bf  Medina A.,  Revoy Ph.,} {\it Groupes de Lie \`a
structure symplectique invariante,} Symplectic geometry, groupoids
and integrable systems, in "S\'eminaire Sud Rhodanien", M.S.R.I.,
New York/Berlin: Springer-Verlag, (1991), 247-266.

[14] {\bf  Medina A.,  Revoy Ph.,} {\it Lattices un Symplectic Lie
Groups,} Journal of Lie Theory, Volume {\bf 17} (2007) 27-39.

[15] {\bf  Onishchik A.L.,  Vinberg E.B.(Eds.),} {\it
 Lie Algebras III. Structure of Lie Groups. and Lie Algebras,}
Encyclopaedia of Mathematical. Sciences, {\bf vol. 41},
Springer-Verlag, Berlin, 1994.

[16] {\bf Tsemo A.,} {\it Automorphismes polynomiaux des
vari\'et\'es affines,}  C. R. Acad. Sci. Paris S\'er. I Math. {\bf
329} (1999), no. 11, 997-1002.
\bigskip

Mohamed Boucetta\\
Facult\'e des Sciences et Techniques \\
BP 549 Marrakech, Morocco.
\\
Email: {\it boucetta@fstg-marrakech.ac.ma }
\bigskip

Alberto Medina\\ Universit\'e Montpellier 2 \\Case Courrier 051,
UMR CNRS 5149\\
Place Eugène Bataillon 34095\\ MONTPELLIER Cedex France\\
Email: {\it medina@math.univ-montp2.fr}

\end{document}